\documentclass[a4paper,11pt]{article} \usepackage{amssymb} 
\usepackage{amsmath} \usepackage{amsfonts} 
\usepackage{amsmath} \usepackage{bm} 
\def\today{}

\def\s01{{\{0,1\}}}
\def\cosec #1 {{\textrm {cosec}\ #1}}

\def\seqq#1#2#3#4{{(#4)_{#1=#2}^{#3}}}
\def\summ#1#2#3{{\sum_{#1=#2}^{#3}}}

\def\cupp#1#2#3{{\cup_{#1=#2}^{#3}}}
\def\capp#1#2#3{{\cap_{#1=#2}^{#3}}}

\def\range #1#2#3{{{#1}={#2},\ldots,{#3}}}
\def\bx{{\bf x} }

\def\by{{\bf y} }
\def\bz{{\bf z} }

\def\b1{{\bf 1} }

\def\del{{\delta}}

\def\alp{{\alpha}}
\def\veps{{\varepsilon}}
\def\bet{{\beta}}
\def\sig{{\sigma}}

\def\gam{{\gamma}}

\def\bN{{\ensuremath{\mathbb N} }}

\def\bQ{{\ensuremath{\mathbb Q} }}
\def\bR{{\ensuremath{\mathbb R} }}

\def\b #1{{\bf #1}}

\def\bu{{\bf u}}
\def\bv{{\bf v}}

\def\proof{{\beginpf}}
\def\endproof{{\eopf}}

\def\vp{{\varphi}}

\def\lin{{\textrm{lin}}}
\def\linbar{{\overline{\textrm{lin}}}}
\def\supp{{\textrm{supp}}}
\newcommand\beq{\begin{equation}}
\newcommand\eeq{\end{equation}}
\newcommand\bdfn{\begin{defn}}
\newcommand\blem{\begin{lemma}}
\newcommand\elem{\end{lemma}}
\newcommand\bcor{\begin{cor}}
\newcommand\ecor{\end{cor}}
\newcommand\bthm{\begin{thm}}
\newcommand\ethm{\end{thm}}
\newcommand\edfn{\end{defn}}
\newcommand\bcas{\begin{cases}}
\newcommand\ecas{\end{cases}}

\def\nm#1{{\left\Vert #1 \right\Vert}}
\def\Nm{{\nm\cdot}}

\def\veps{{\varepsilon}}

\newcommand\sprod[2]{\langle{#1},{#2}\rangle}

\def\casif#1!{& \textrm{#1}\\}
\newif\ifrough
\roughfalse

\def\reff#1{{\ref {#1}}}
\def\namedlabel#1!{{\label{#1}}}

\def\refeq#1{{(\reff{#1})}}
\def\up#1{^{(#1)}}

\def\reflem#1{{Lemma \reff{#1}}}

\def\refthm#1{{Theorem \reff{#1}}}

\def\benum{\begin{enumerate}}
\def\eenum{\end{enumerate}}

\newtheorem{thm}{Theorem}[section]
\newtheorem{defn}[thm]{Definition}

\newtheorem{cor}[thm]{Corollary}

\newtheorem{lemma}[thm]{Lemma}

\def\ap1{{\square}}

\newcounter{smallromans}

\newcounter{smallarabics}

\newcounter{smallalphs}

\newcommand\bproof{\begin{proof}}
\newcommand\eproof{\end{proof}}
\newcommand{\beginpf}{\smallskip\textbf{Proof. }}
\newcommand{\eopf}{\hfill $\Box$}

\roughfalse

\begin{document} 
\title{ Banach spaces with no proximinal subspaces of codimension 2 }  
\author{C.J.Read\\ Department of Pure Mathematics, University of Leeds, LS2 
9JT, UK\\ solocavediver@gmail.com; http://solocavediver.com/maths}
\date{\today } 
\maketitle 

\abstract{
	The classical theorem of Bishop-Phelps asserts that, for a Banach space $X$,
the norm-achieving functionals in $X^*$ are dense in $X^*$. B\'ela Bollob\'as's extension of the theorem gives a quantitative description of just how dense the norm-achieving functionals have to be: if $(x,\vp)\in X\times X^*$ with $\nm x=\nm {\vp}=1$ and $|1-\vp(x)|<\veps^2/4$ then there are 
$(x',\vp')\in   X\times X^*$ with $\nm {x'}=\nm {\vp'}=1$, $\nm{x-x'}\vee \nm{\vp-\vp'}<\veps$  and 
$\vp'(x')=1$.

This means that there are always ``proximinal'' hyperplanes $H\subset X$ (a nonempty subset $E$ of a metric space is said to be ``proximinal'' if, for $x\notin E$, the distance $d(x,E)$ is always achieved - there is always an $e\in E$ with $d(x,E)=d(x,e)$); for if $H=\ker\vp$ ($\vp\in X^*$) then it is easy to see that $H$ is proximinal if and only if $\vp$ is norm-achieving. Indeed the set of proximinal hyperplanes $H$ is, in the appropriate sense, dense in the set of all closed hyperplanes $H\subset X$.

Quite a long time ago [Problem 2.1 in his monograph ``The Theory of Best approximation and Functional Analysis" Regional Conference series in Applied Mathematics, SIAM, 1974], Ivan Singer
 asked if this result generalized to closed subspaces of finite codimension -  if every Banach space has a proximinal subspace of codimension 2, for example. In this paper I will show that there is a Banach space $X$ such that $X$ has no proximinal subspace of finite codimension $n\ge 2$.  So we have a converse to Bishop-Phelps-Bollob\'as: a dense set of proximinal hyperplanes can always be found, but proximinal subspaces of larger, finite codimension need not be.

}

\vfill\eject

\section{Introduction.}\namedlabel 1!
I'm grateful to David Blecher for awakening me to the joys of proximinality in the context of operator algebras (norm-closed subalgebras of $B(H)$), and to Gilles Godefroy for alerting me to this particular problem.  

The original Bishop-Phelps theorem is \cite{BP}, and Bollob\'as' improved version of the theorem is \cite{BB}. 
The place where the problem solved in this paper was originally posed is in Ivan Singer \cite{S}. Gilles Godefroy's exhaustive survey article on isometric preduals  in Banach spaces, which discusses this problem among many others, is \cite{GG}. Our work with David Blecher involving proximinality of ideals in operator algebras is \cite{BR}. This is a successful attempt to generalize, to a noncommutative setting, the classical Glicksberg peak set theorem in uniform algebras (Theorem 12.7 in Gamelin \cite{G}).

All the Banach spaces in this paper are over the real field. At risk of stating the obvious, a proximinal subset is necessarily closed; so we lose no generality later on by assuming that a (hypothetical) proximinal subspace of finite codimension is the intersection of the kernels of finitely many continuous linear functionals.

Let $c_{00}(\bQ)$ denote the terminating sequences with rational coefficients (a much-loved countable set), and let $\seqq k1\infty{\bu_k}$ be a sequence of elements of $c_{00}(\bQ)$ which lists every element infinitely many times. For $\bx\in c_{00}(\bQ)$, write $u^{-1}\{\bx\}$ for the infinite set $\{k\in\bN:\bu_k=\bx\}$.

Let $\seqq k1\infty{a_k}$ be a strictly increasing sequence of positive integers. We impose a growth condition: if $\bu_k\ne 0$, we demand that 
\beq\namedlabel 1-1!
a_k>\max\ \supp\ \bu_k, \textrm{and}\ a_k\ge\nm{\bu_k}_1,
\eeq 
where $\supp\ \bu$ denotes the (finite) support of $\bu\in c_{00}(\bQ)$, and $\nm{\bu}_1$ denotes the $l^1$ norm. 

For $E\subset \bN$ we write $A_E$ for the set $\{a_k:k\in E\}$; for $\bx\in c_{00}(\bQ)$ we write 
$A_\bx$ for $A_{\bu^{-1}\{\bx\}}$. $A_\bx$  is an infinite set, and in view of \refeq{1-1}, for each $\bx\in c_{00}(\bQ)\setminus\{0\}$ we have 
\beq\namedlabel 1-2!
\min A_\bx>\max\ \supp\ \bx,\ \ \ \min A_\bx\ge\nm{\bx}_1.
\eeq

Given sequences $(\bu_k)$, $(a_k)$ as described above, we define a new norm $\Nm$ on $c_0$ as follows:
\beq\namedlabel 2-1!
\nm\bx={\nm\bx}_0+\summ k1\infty 2^{-a_k^2}|\sprod\bx{\bu_k-e_{a_k}}|.
\eeq
Here ${\nm\bx}_0=\sup_n|x_n|$ is the usual norm on $c_0$; $(e_j)$ are the unit vectors; and the duality $\sprod\bx{\bu_k-e_{a_k}}$ is the $\sprod{c_0}{l^1}$ duality. Now in view of \refeq{1-1}, we have $\nm{\bu_k-e_{a_k}}_1=1+\nm{\bu_k}_1\le 1+a_k$, so $\summ k1\infty 2^{-a_k^2}{\nm{\bu_k-e_{a_k}}}_1\le \summ k1\infty 2^{-a_k^2}(1+a_k)\le \summ n1\infty(1+n)\cdot 2^{-n^2}<2$. Accordingly, we have 
\beq\namedlabel 2-2!
\nm{\bx}_0\le\nm\bx\le 3\nm{\bx}_0
\eeq
for all $\bx\in c_0$. For our main theorem in this paper, we shall show:

\bthm\namedlabel t2!
The Banach space $(c_0,\Nm)$ has no proximinal subspace $H$ of finite codimension $n\ge 2$.
\ethm
\def\gat{{G\^ateaux}}

\section{\gat\ derivatives}\namedlabel 2!

Recall that if $X$ is a real vectorspace, $\bu,\bx\in X$ and $f:X\to\bR$, then the \gat\ derivative (of $f$, at $\bx$, in direction $\bu$) is defined as
\beq\namedlabel 2-3!
df(\bx;\bu)=\lim_{h\to 0}\frac{f(\bx+h\bu)-f(\bx)}h,
\eeq
when that limit exists. We will make use of the one-sided forms of this derivative:
\beq\namedlabel 2-4!
df_+(\bx;\bu)=\lim_{h\to 0+}\frac{f(\bx+h\bu)-f(\bx)}h,
\eeq
and
\beq\namedlabel 2-5!
df_-(\bx;\bu)=\lim_{h\to 0-}\frac{f(\bx+h\bu)-f(\bx)}h.
\eeq
Obviously $df_-(\bx;\bu)=-df_+(\bx;-\bu)$ for all $f,\bx$ and $\bu$ such that either derivative exists.
Of particular interest to us is when $X=c_0$ and $f(\bx)=\nm \bx$ as defined in \refeq{2-1} (the ``usual'' norm on $c_0$ will always be referred to as $\Nm_0$ in this paper). The derivatives 
$d_\pm  f(\bx;\bu)$ for this function $f$ will be written $d_\pm\nm{\bx;\bu}$. 
Now it is a fact that the derivative $d_\pm\nm{\bx;\bu}$ exists everywhere. To see this, let us prove some small lemmas.

\blem\namedlabel 3.1!
If ${\nm\bx}_0$ denotes the $c_0$-norm, the derivative $d_+\nm{\bx;\bu}_0$ exists at all points $\bx,\bu\in c_0$. In fact, if $\bx=0$ then the derivative is $\nm\bu_0$; whereas if $\bx\ne 0$, we may write $E_+=\{n\in\bN:|x_n|=\nm\bx_0, u_nx_n>0\}$ and  $E_-=\{n\in\bN:|x_n|=\nm\bx_0, u_nx_n\le 0\}$, and we have 
\beq\namedlabel 3-1!
d_+\nm{\bx;\bu}_0=\bcas
\max\{|u_n|:n\in E_+\}, &\textrm{if } E_+\ne\emptyset;\cr
-\min\{|u_n|:n\in E_-\}, &\textrm{if } E_+=\emptyset, E_-\ne\emptyset.\cr
\ecas
\eeq
\elem\proof This is an easy calculation which we omit (note that $E_+$ and $E_-$ cannot both be empty!). 

\blem\namedlabel 4.1!
Let $X$ be a Banach space and $\vp\in X^*$. Then the \gat\ derivative of $f(\bx)=|\vp(\bx)|$ exists at all points $(\bx;\bu)\in X\times X$. We have 
$$
d_+f(\bx;\bu)=f(\bu)\sig(\vp(\bu)\vp(\bx)),
$$ 
where the sign
\beq\namedlabel 4-1!
\sig(t)=\bcas
+1,&\textrm{if } t\ge 0;\cr
-1,&\textrm{if } t<0.
\ecas
\eeq\elem
\proof This is an even simpler calculation, which we also omit.

\def\lip{{\textrm{Lip}}}

\bdfn\namedlabel 4.2!
For a real normed space $X$ and a function $f:X\to\bR$, define the Lipschitz constant\edfn
\beq\namedlabel 4-2!
\lip_1f=\sup\{\frac{|f(\bx)-f(\by)|}{\nm{\bx-\by}}:\bx,\by\in X, \bx\ne\by\}.
\eeq

\blem\namedlabel 4.3!
Let $X$ be a real normed space, and $\seqq n0\infty{f_n}$ a sequence of functions from $X$ to $\bR$, such that $d_+f_n(\bx;\bu)$ exists at each $(\bx;\bu)\in X\times X$. Suppose $\summ n0\infty\lip_1 f_n<\infty$, and $\summ n0\infty f_n(0)$ converges. Then the function $f=\summ n0\infty f_n$ exists everywhere on $X$, and $d_+f=\summ n0\infty$ $d_+f_n$ exists everywhere also.
\elem
\proof The sum $\summ n0\infty f_n(\bx)$ converges because $\summ n0\infty f_n(0)$ converges, and $|f_n(\bx)-f_n(0)|\le\nm\bx\cdot\lip_1f_n$ so $\summ n0\infty f_n(\bx)-f_n(0)$ converges also. Since $|d_+f_n(\bx;\bu)|\le\nm\bu\cdot\lip_1f_n$, we find that the sum $\summ n0\infty d_+f_n(\bx;\bu)$
converges; we claim the sum is $d_+f(\bx;\bu)$. For given $\bx,\bu\ne 0,$ and $\veps>0$, we can choose $N$ so large that $\nm u\cdot\summ n{N+1}\infty\lip_1 f_n<\veps/3$, so 
\beq\namedlabel 5-1!
|\summ n{N+1}\infty d_+f_n(\bx;\bu)|<\veps/3\eeq
and for every $h>0$,
\beq\namedlabel 5-2!
|\summ n{N+1}\infty(f_n(\bx+h\bu)-f_n(\bx))/h|\le 
\nm\bu\cdot\summ n{N+1}\infty\lip_1f_n<\veps/3
\eeq
also. As $h\to 0+$, we know 
$(\summ n1Nf_n(\bx+h\bu)-\summ n1Nf_n(\bx))/h \to \summ n1Nd_+f_n(\bx;\bu)$, so we can choose $\del>0$ such that whenever $0<h<\del$, we have 
\beq\namedlabel 5-3!
|(\summ n1Nf_n(\bx+h\bu)-\summ n1Nf_n(\bx))/h - \summ n1Nd_+f_n(\bx;\bu)|<\veps/3.
\eeq 
Adding up \refeq{5-1}, \refeq{5-2} and \refeq{5-3}, we find that whenever $0<h<\del$, we have 
$|(f(\bx+h\bu)-f(\bx))/h-\summ n1\infty d_+f_n(\bx;\bu)|<\veps$. This completes the proof.\endproof

\bcor\namedlabel c5!
The new norm $\Nm$ on $c_0$ has a one-sided derivative $d_+\nm{\bx;\bu}$ everywhere. 
Furthermore,
\beq\namedlabel 5-5!
d_+\nm{\bx;\bu}=d_+\nm{\bx;\bu}_0+\summ k1\infty 2^{-a_k^2}\sig_k|\sprod\bu{\bu_k-e_{a_k}}|,
\eeq
where $\sig_k=\sig_k(\bx;\bu)=\sig(\sprod{\bu}{\bu_k-e_{a_k}}\sprod{\bx}{\bu_k-e_{a_k}})$,
and the function $\sig$ is as in \refeq{4-1}.
\ecor

\proof If we write $f_0(\bx)=\nm\bx_0$ and $f_k(\bx)=2^{-a_k^2}|\sprod\bx{\bu_k-e_{a_k}}|$, then the Lipschitz constants for $f_k$ are 1 (if $k=0$) or 
$2^{-a_k^2}\nm{\bu_k-e_{a_k}}_1\le(1+a_k)\cdot 2^{-a_k^2}$ for $k>0$. Accordingly $\summ k0\infty\lip_1f_k<\infty$, and the derivatives $d_+f_k$ are given by \reflem{3.1} and \reflem{4.1}. We have $\nm\bx=\summ k0\infty f_k(\bx)$ so $d_+\nm{\bx;\bu}=\summ n0\infty d_+f_n(\bx;\bu)$
  by \reflem{4.3}. This sum works out to expression \refeq{5-5}.\endproof

The key link between \gat\ derivatives and proximinality is as follows:

\blem\namedlabel 6.1!
Suppose $(X,\Nm)$ is a Banach space, $H\subset X$ a subspace, and suppose that for some $\bx\in X\setminus H$, and $\bv\in H$, the \gat\ derivatives $d_\pm\nm{\bx;\bv}$ both exist, are nonzero, and have the same sign. Then $\nm x\ne\inf\{\nm \by:\by\in\bx+H\}$. $\bx$ is not a closest point to zero in the coset $\bx+H$.
\elem
\proof We may consider $\by=\bx+h\bv$ for small nonzero $h\in\bR$. Depending on the sign of $h$, the norm $\nm \by$ is roughly $\nm\bx+h\cdot d_\pm \nm{\bx;\bv}$. But the signs of $d_\pm\nm{\bx;\bv}$ are the same, so if $h$ is chosen correctly, we get $\nm \by<\nm\bx$.
\endproof

\bcor\namedlabel 6.2!
Suppose $H\subset X$ as in \reflem{6.1}, and there is an $\bx\in X\setminus H$ such that for every $\bz\in H$, there is a $\bv\in H$ such that the \gat\ derivatives $d_\pm\nm{\bx+\bz;\bv}$ exist, are nonzero, and have the same sign. Then $H$ is not proximinal in $X$.
\ecor
\proof For in this case, there is no element $\bx+\bz\in\bx+H$ which achieves the minimum distance from that coset to zero. Equivalently, there is no element $\bz\in H$ which achieves the minimum distance from $H$ to $-\bx$. $H$ is not proximinal.
\endproof

\section{Approximate linearity of $d_\pm$}\namedlabel 3!

It is a feature of the \gat\ derivative $df(\bx;\bu)$ that it does not have to be linear in $\bu$. This is of course also true of the single-sided derivatives $df_\pm$. So, in this section we develope a result asserting ``approximate linearity'' of $d_\pm\nm{\bx;\bv}$ for $\bx,\bv\in c_0$.

\bdfn\namedlabel 7.1!
Let $f:c_0\to\bR$ be such that $d_+f(\bx;\bv)$ exists for all $\bx,\bv\in c_0$. Let $\bx\in c_0$, and let $\gam\in l^1$ be such that the support $E=\{i:\sprod {e_i}\gam\ne 0\}$ is infinite. We shall say $d_+f(\bx)$ is approximately linear on $E$ (and approximately equal to $\gam$) if there is an ``error  sequence'' $(\veps_i)_{i\in E}$ with $\veps_i>0$, $\veps_i\to 0$, such that for all $\bv\in c_0$ with $\supp\ \bv\subset E$, we have 
\edfn 
\beq\namedlabel 7-1!
|d_+f(\bx;\bv)-\sprod\bv\gam|\le\sum_{i\in E}\veps_i|v_i\gam_i|.
\eeq
Note that if $\bv$ is chosen so that $|\sprod\bv\gam|> \sum\veps_i|v_i\gam_i|$, then \refeq{7-1} implies that $d_+f(\bx;\bv)$ is nonzero, and has the same sign as $d_-f(\bx;\bv)=-d_+f(\bx;-\bv)$.

\blem\namedlabel 7.2!
Let $\bx\in c_0$ be given, and $\bz_1,\ldots,\bz_m\in c_{00}$ such that $\sprod\bx{\bz_j}\ne 0$ for any $\range j1m$. Let $A_{\bz_i}=A_{u^{-1}\{\bz_i\}}$ as in \S\reff 1, and let $A=\cupp i1mA_{\bz_i}$. 
Then $d_+\nm\bx$ is approximately linear on a cofinite subset $A_0\subset A$, the derivative being approximately equal to $\gam=\summ i1\infty\gam_i e_i^*$, where
\elem
\beq\namedlabel 7-2!
\gam_i=
\bcas
-2^{-a_k^2}\sig(\sprod\bx{\bz_j})&\textrm{if } i=a_k\in A_0\cap A_{\bz_j}
\cr
0,&\textrm{otherwise}.
\ecas
\eeq
The error sequence $(\veps_i)_{i\in A_0}$ can be taken to be
\beq\namedlabel 7-3!
\veps_i=2^{a_k^2}\cdot\summ l{k+1}\infty 2^{-a_l^2},\ i=a_k\in A_0.
\eeq
\proof Let $\bv$ be any vector supported on $A$. The error $\del=d_+\nm{\bx;\bv}-\sprod\bv\gam$
is given by \refeq{5-5}\ and \refeq{7-2}; we have 
$$\del=d_+\nm{\bx;\bv}_0+\sum_{k\in\bN\setminus\cupp j1m u^{-1}\{\bz_j\}}2^{-a_k^2}\sig_k|\sprod\bv{\bu_k-e_{a_k}}|
$$
\beq\namedlabel 8-1!
+\summ j1m\sum_{k\in u^{-1}\{\bz_j\}}
2^{-a_k^2}\sig_k|\sprod\bv{\bu_k-e_{a_k}}|+v_{a_k}\sig(\sprod\bx{\bz_j})
\eeq
where $\sig_k=\sig(\sprod\bv{\bu_k-e_{a_k}}\sprod\bx{\bu_k-e_{a_k}})$.

Now by \reflem{3.1}, the derivative $d_+\nm{\bx;\bv}_0$ is zero unless $v_i\ne 0$ for some $i\in E=\{n:|x_n|=\nm{\bx}_0\}$. This set $E$ is finite, and $\bv$ will be supported on $A_0$; so if we choose our cofinite set $A_0\subset A$ so that $A_0\cap E=\emptyset$, we have 
 $d_+\nm{\bx;\bv}_0=0$. If we also ensure that $A_0\cap\supp\ \bz_j=\emptyset$ for each $\range j1m$, we find that when $\bv$ is supported on $A_0$, and $k\in u^{-1}\{\bz_j\}$, we have 
\beq\namedlabel 8-1b!
\sprod\bv{\bu_k-e_{a_k}}=
\sprod\bv{\bz_j-e_{a_k}}=-\sprod\bv{e_{a_k}}=-v_{a_k}. 
\eeq
So if we choose $A_0$ so that $A_0\cap(E\cupp j1m\supp\ \bz_j)=\emptyset$, the expression
\refeq{8-1} simplifies somewhat to 
$$\del=\sum_{k\in\bN\setminus\cupp j1m u^{-1}\{\bz_j\}}
2^{-a_k^2}\sig_k|\sprod\bv{\bu_k-e_{a_k}}|
$$
\beq\namedlabel 8-2!
+\summ j1m\sum_{k\in u^{-1}\{\bz_j\}}(\sig_k|v_{a_k}|+v_{a_k}\sig(\sprod\bx{\bz_j}));
\eeq
and $\sig_k$ itself simplifies to $\sig_k=\sig(-\sprod\bx{\bz_j-e_{a_k}}\cdot v_{a_k})$. 
Now $F_j=\{k\in\bN:|\sprod\bx{e_{a_k}}|\ge |\sprod\bx{\bz_j}|\}$ is a finite set; we may thus also assume that $A_0$ does not meet any $F_j$. In that case, $\sprod \bx{\bz_j-e_{a_k}}$ is nonzero and has the same sign as $\sprod\bx{\bz_j}$, so $\sig_k=-v_{a_k}\sig(\sprod\bx{\bz_j})$. So the second term in \refeq{8-2}\ disappears, and we have 
\beq\namedlabel 9-1!
\del=\sum_{k\in\bN\setminus\cupp j1m u^{-1}\{\bz_j\}}2^{-a_k^2}\sig_k|\sprod\bv{\bu_k-e_{a_k}}|.
\eeq
Even better, $\bv$ is supported on $A=a\cdot\cupp j1m u^{-1}\{\bz_j\}$, so all terms $\sprod\bv{e_{a_k}}$ are zero in \refeq{9-1}, and we have 
$$
\del=\sum_{l\in\bN\setminus\cupp j1m u^{-1}\{\bz_j\}}2^{-a_l^2}\sig_l|\sprod\bv{\bu_l}|;
$$
\beq\namedlabel 9-1b!
|\del|\le\sum_{i\in A_0}\sum_{l\in\bN\setminus\cupp j1m u^{-1}\{\bz_j\}}2^{-a_l^2}
|v_i|\cdot |\sprod{e_i}{\bu_l}|.
\eeq
Now in every case when $i\in A_0$ we have $i=a_k$ for some $k\in u^{-1}\{\bz_j\}$, $\range j1m$. If $\sprod {e_{a_k}}{\bu_l}\ne 0$ then the support $\supp\ \bu_k$ is not contained in $[0,a_k)$.
But if  $l\le k$ then the support of $\bu_l$ is contained in $[0,a_k)$ by \refeq{1-1}. So $k<l$ in every case when $\sprod {e_{a_k}}{\bu_l}\ne 0$. Accordingly, 
$$
|\del|\le\sum_{i=a_k\in A_0}\sum_{l>k}2^{-a_l^2}|v_i|\cdot |\sprod{e_i}{\bu_l}|
\le \sum_{i=a_k\in A_0}\sum_{l>k}2^{-a_l^2}\cdot a_l\cdot |v_i|
$$
because $\nm{\bu_l}_1\le a_l$ by \refeq{1-1} again. Now for $i=a_k\in A_0$, we have 
$|\gam_i|=2^{-a_k^2}$ by \refeq{7-2}; so writing $\veps_i=2^{a_k^2}\cdot\summ l{k+1}\infty 2^{-a_l^2}$ as in \refeq{7-3}, we have $\veps_i\to 0$ and 
$$
|\del|=|d_+\nm{\bx;\bv}-\sprod\bv\gam|\le\sum_{i\in A_0}\veps_i|v_i\gam_i|
$$
exactly as in \refeq{7-1}. So $d_+\nm\bx$ is approximately linear on a cofinite subset $A_0\subset A$, with the derivative $\gam\in l^1$ given by \refeq{7-2}, and the error sequence $(\veps_i)$ given by \refeq{7-3}. \endproof

\section{Using the Hahn-Banach Theorem.}\namedlabel 4!

\blem\namedlabel 10!
Let $H\subset c_0$ be a closed subspace of finite codimension, say $H=\capp i1N\ker\vp_i$, where each $\vp_i\in l^1$. Let $\bx\notin H$ be an element of minimum norm in the coset $\bx+H$, and let $\bz_1,\ldots, \bz_m\in c_{00}$ be such that $\sprod\bx{\bz_j}\ne 0$ for any $\range j1m$. Let $A_0\subset A=\cupp j1mA_{\bz_j}$ be a cofinite subset satisfying the conditions of \reflem{7.2}, and let $\gam\in l^1$ be the approximate derivative as in \reflem{7.2}, $(\veps_i)_{i\in A_0}$ the error sequence as in \refeq{7-3}. Then there is a $\vp\in\lin\{\vp_j:i\le j\le m\}$ such that for every $i\in A_0$,
we have 
\beq\namedlabel 10-0!
|\sprod{e_i}\vp-\gam_i|\le\veps_i|\gam_i|.
\eeq
\elem 
%$(e_i^*)_{i\in\bN}$ denote the unit vector basis of $l^1$. 
\proof We consider the weak-* topology on $l^1$ with respect to its usual predual, $c_0$.  The set $G=\{\vp\in l^1: |\sprod{e_i}\vp-\gam_i|\le\veps_i|\gam_i|$
(all $i\in A_0),$ and $\sprod {e_i}\vp=0$ (all $i\notin A_0)\}$ is a weak-* compact convex set. The set
 $\Phi=\lin\{\vp_i,\range i1N\}+\linbar\{e_j:j\in\bN\setminus A_0\}\subset l^1$ is a weak-* closed subspace, because it is $\{\vp\in l^1:\vp(\bu)=0$ for every $\bu\in c_0$ supported on $A_0$, such that $\vp_i(\bu)=0$ ($\range i1N)\}$. 

If $\Phi\cap G\ne\emptyset$, then the assertion of the Lemma is satisfied. If $\Phi\cap G=\emptyset$, then the Hahn-Banach Separation Lemma tells us that there is a weak-* continuous $\bv \in l^\infty$ separating them; of course the weak-* continuity means that $\bv\in c_0$.  We may assume $\sprod\vp\bv=0$ for $\vp\in\Phi$, but $\sprod\vp\bv\ge 1$  whenever $\vp\in G$. Since $\bv$ annihilates $\Phi$, the support of  
$\bv$ is contained in $A_0$. By approximate linearity of $d_+\nm\bx$, from \refeq{7-1} we have 
\beq\namedlabel 11-2!
|d_+\nm{\bx;\bv}-\sprod\bv\gam|\le\sum_{i\in A_0}\veps_i|v_i\gam_i|;
\eeq
and the same is true with $d_+$ replaced by $d_-$. We cannot have $d_+\nm{\bx;\bv}$ and 
$d_-\nm{\bx;\bv}$ the same sign, or \reflem{6.1} would tell us $\bx$ does not have minimum norm in the coset $\bx+H$. So, as observed after \refeq{7-1}, we must have 
\beq\namedlabel 11-3!
|\sprod\bv\gam|\le\sum_{i\in A_0}\veps_i|v_i\gam_i|.
\eeq
Let us write $\eta=\sprod\bv\gam/\sum_{i\in A_0}\veps_i|v_i\gam_i|\ \in[-1,1]$ (noting that the denominator cannot be zero since $\veps_i,\gam_i$ are never zero for $i\in A_0$, and $\bv\ne 0$ is supported on $A_0$). Define a new $\vp\in l^1$ by
\beq\namedlabel 12-1!
\sprod{e_i}\vp=\bcas
\gam_i(1-\eta\veps_i\sig(v_i\gam_i)), &$if$\ i\in A_0;\cr
0,&$otherwise.$
\ecas
\eeq
We then have $|\sprod{e_i}\vp-\gam_i|\le\veps_i|\gam_i|$ ($i\in A_0$), so $\vp\in G$, yet
\beq\namedlabel 12-2!
\sprod\bv\vp=\sprod\bv\gam-\eta\cdot\sum_{i\in A_0}\veps_iv_i\gam_i\sig(v_i\gam_i)=0.
\eeq
This contradicts the Hahn-Banach separation of $\bv$, which asserts that for such $\vp$ we should have $\sprod\bv\vp\ge 1$. Thus the Lemma is proved.\endproof

Let us now begin to use our information to investigate proximinal subspaces.
If $i=a_l$ for some $l\in\bN$, we shall write $\alp_i=2^{-a_l^2}$.

\bthm\namedlabel 12!
Let $H\subset (c_0,\Nm)$ be a proximinal subspace of finite codimension,
say $H=\capp i1N\ker\vp_i$, where $\vp_i\in l^1$. 
Let $\bz_j\in c_{00}$ ($\range j1m$), and $A=\cupp i1mA_{\bz_i}$.
Write $\Phi=\lin\{\vp_i:\range i1N\}$, and 
let
\beq
\Phi_0=\{\vp\in\Phi:\sup\{\alp_i^{-1}|\sprod{e_i}\vp|:i\in A\}<\infty\}.
\eeq
 Let $\theta_0:\Phi_0\to l^\infty(A)$ be the linear map such that
\beq\namedlabel 12-5!
(\theta_0\vp)_i=\alp_i^{-1}\sprod{e_i}\vp \ (i\in A);
\eeq
and let $q:l^\infty(A)\to l^\infty(A)/c_0(A)$ be the quotient map.
Write $\theta=q\theta_0$. Let $\bx\in H$ be an element such that $\nm\bx$ is minimal in the coset $\bx+H$, and suppose $\sprod\bx{\bz_i}\ne 0$ for any $\range i1m$. Then the image $\theta\Phi_0$ includes the vector 
$\sig_\bx+c_0(A)\in l^\infty(A)/c_0(A)$, where 
\ethm
\beq\namedlabel 13-1!
(\sig_\bx)_i=\sig(\sprod\bx{\bz_j}) \textrm{ if } i\in A_{\bz_j}, \range j1m.
\eeq
\proof
By \reflem{10}, there is a $\vp\in\Phi$ such that for all but finitely many $i\in A$, we have 
$|\sprod{e_i}\vp-\gam_i|\le \veps_i|\gam_i|$, where 
for $i=a_l\in A_{\bz_j}$ we define $\gam_i=-2^{-a_l^2}\sig(\sprod\bx{\bz_j})$ $=-\alp_i\sig(\sprod\bx{\bz_j})$. Since $\veps_i\to 0$ it is clear that $\sup\{\alp_i^{-1}|\sprod{e_i}\vp|\}<\infty$, so $\vp\in\Phi_0$, and the image $\theta\vp$ is the vector 
$-\sig_\bx+c_0(A)$, because for $i=a_l\in A_{\bz_j}$ we have 
$(\theta_0\vp)_i=\alp_i^{-1}\sprod{e_i}\vp\in -\sig(\sprod\bx{\bz_j})+[-\veps_i,\veps_i]$.
So, $\theta(-\vp)=\sig_\bx+c_0(A)$.\endproof 

\section{Proof of \refthm{t2}}\namedlabel 5!

Suppose towards a contradiction that $H\subset(c_0,\Nm)$ is a proximinal subspace of finite codimension $N\ge 2$. Any proximinal subspace must be closed, so let us say $H=\capp i1N\ker\vp_i$, where the $\vp_i\in l^1$ are linearly independent. For $\range r0{N+1}$, let us write
$\bet_r=r\pi/(2N+2)$, and for $\range r1{N+1}$ let us pick $\bx\up r\in c_0$ such that $\sprod{\bx\up r}{\vp_1}=\cos\bet_r$, and $\sprod{\bx\up r}{\vp_2}=\sin\bet_r$. Perturbing each $\bx\up r$ by an element of $H$ as necessary, we can assume that each $\nm{\bx\up r}$ is minimal in the coset 
$\bx\up r+H$. Writing $\zeta_r=(\bet_r+\bet_{r-1})/2$ ($\range r1{N+1}$), we define the linear functional $\psi_r=\sin\zeta_r\cdot\vp_1-\cos\zeta_r\cdot\vp_2$, so
$$
\sprod{\bx\up r}{\psi_s}=\sin\zeta_s\cos\bet_r-\cos\zeta_s\sin\bet_r=\sin(\zeta_s-\bet_r);
$$
thus $\sprod{\bx\up r}{\psi_s}>0$ if $s>r$, but $\sprod{\bx\up r}{\psi_s}<0$ if $s\le r$.

Pick a finite sequence $\seqq r1{N+1}{\bz_r}\subset c_{00}$ with $\nm{\bz_r-\psi_r}_1$ sufficiently 
small, and we will also find that $\sprod{\bx\up r}{\bz_s}>0$ if $s>r$, but $\sprod{\bx\up 
r}{\bz_s}<0$ if $s\le r$. We find that the sequence $(\sig(\sprod{\bx\up r}{\bz_s}))_{s=1}^{N+1}$ 
$\in \bR^{N+1}$ is the vector $\by_r=(-1,-1,\ldots,-1,1,1,\ldots 1)$, where there are $r$ entries $-1$ followed 
by $N+1-r$ entries $+1$. It is a fact that the $\by_r$ span $\bR^{N+1}$ - they are linearly independent.

We can apply \refthm{12}\ with the sequence $\bz_1,\ldots,\bz_{N+1}$, and $\bx$ can be any of the 
vectors $\bx\up 1, \ldots, \bx\up {N+1}$. The map $\theta$ is the same for each $\bx\up r$ (because 
the sequence $\alp_i$ doesn't change, only the signs $\sig(\sprod{\bx\up r}{\bz_s})$). Writing $A=\cupp j1{N+1}A_{\bz_j}$, we find that the image 
$\theta\Phi_0$ must contain, for each $\range r1{N+1}$, the vector $\sig_{\bx\up r}+c_0(A)$ with 
\beq
(\sig_{\bx\up r})_i=\sig(\sprod{\bx\up r}{\bz_j})=\sprod{\by_r}{e_j}\textrm{ for all }i\in A_{\bz_j}.
\eeq
(where here $\seqq j1{N+1}{e_j}$ denote the unit vector basis of $\bR^{N+1}$).
Because the vectors $\by_r$ are independent, the dimension of $\theta\Phi_0$ must be at least 
$N+1$. However $\Phi_0\subset\Phi$, and $\dim\Phi=N$. This contradiction implies that $H$ is not proximinal. \endproof

\end{document}